\documentclass[12pt]{article}

\usepackage{epsf,epsfig,amsmath,amsfonts,amsthm}
\usepackage{graphicx}
\usepackage{color}
\usepackage{accents}
\usepackage{mathrsfs}
\usepackage{tikz}
\usepackage{hyperref}

\newtheorem{thm}{Theorem}[section]

\newtheorem{cor}[thm]{Corollary}

\newtheorem{lem}[thm]{Lemma}

\newtheorem{ass}{Assumption}[section]

\numberwithin{equation}{section} 

\usetikzlibrary{patterns.meta}
\usetikzlibrary{calc}
\usetikzlibrary{arrows.meta}
\usetikzlibrary{hobby}
\usepackage{pgfplots}
\usepackage{pifont}
\usepackage{caption}
\usepgfplotslibrary{fillbetween}

\newcommand{\be}{\begin{equation*}}
\newcommand{\ee}{\end{equation*}}
\newcommand{\ben}{\begin{equation}}
\newcommand{\een}{\end{equation}}

\newcommand{\EXP}{\operatorname{\textsf{\upshape E}}}
\newcommand{\PR}{\operatorname{\textsf{\upshape P}}}

\newcommand{\di}{\mathrm{d}}

\newcommand{\inv}{{\rm inv}}

\newcommand{\bbr}{\mathbb{R}}

\newcommand{\fcal}{\mathcal{F}}

\newcommand{\lcal}{\mathcal{L}}

\newcommand{\tcal}{\mathcal{T}}

\newcommand{\qcal}{\mathcal{Q}}

\newcommand{\cscr}{\mathscr C}
\newcommand{\dscr}{\mathscr D}
\newcommand{\lscr}{\mathscr L}

\newcommand{\sscr}{\mathscr S}
\newcommand{\wscr}{\mathscr W}

\usepackage{color}

\begin{document}

\title{On the applicability of the maximality principle \\ 
in optimal stopping: an example \\
involving a diffusion and its running maximum}

\author{
{\sc Neofytos Rodosthenous\footnote{Department of
Mathematics, University College London, Gower Street,
London WC1E 6BT, UK,
\texttt{n.rodosthenous@ucl.ac.uk} }
\ and Mihail Zervos\footnote{Department of Mathematics,
London School of Economics, Houghton Street, London
WC2A 2AE, UK, \texttt{mihalis.zervos@gmail.com}}}}

\maketitle

\begin{abstract} \noindent
The maximality principle has proved to be a valuable tool
for identifying the optimal stopping boundaries in
optimal stopping problems involving one-dimensional
time-homogeneous diffusions and their running maximum
processes.
The principle characterises the optimal stopping boundary
as the maximal solution to a first-order nonlinear ODE
that remains strictly below the diagonal in $\bbr^2$ or
a ``zero-level curve''.
In this paper, we construct a suitably tailored optimal
stopping problem to which the maximality principle
is not applicable. 
\\\\
\noindent
{\em MSC2010 subject classification\/}: 60G40, 60H30,
49J10, 49K10, 93E20
\end{abstract}

\section{Introduction}

In several optimal stopping problems involving a
one-dimensional time-homogeneous diffusion $X$ and
its running maximum process $S$, the solution to
a suitable first-order nonlinear ODE determines
(part of) the boundary separating the so-called
stopping and waiting regions.
In some cases, the structure of the problem yields
a natural boundary condition or asymptotic behaviour
that determines the solution to this ODE corresponding
to the optimal stopping boundary.
This is the case for the optimal stopping problems
studied by
Dubins, Shepp and Shiryaev~\cite{DSS94},
Guo and Shepp~\cite{GS01},
Gapeev~\cite{G07},
Ott~\cite{O13, O14}, and
Rodosthenous and Zhang~\cite{RZ20}.
However, in other cases, such conditions do not arise
naturally.
In a seminal paper, Peskir~\cite{P98} provided a
systematic way of selecting the relevant solution
among the continuum of solutions to this ODE.
He proved that the optimal free-boundary
function is given by the unique solution
that satisfies what he termed the {\em maximality
principle\/}. 
In particular, he showed that, for problems in which the
maximality principle holds, the optimal stopping boundary
can be derived and completely characterised solely by
an analysis of its associated ODE.

Since then, the maximality principle has been
used or observed to hold true in
numerous research contributions.
For instance, see
Graversen and Peskir~\cite{GP98a, GP98b},
Peskir~\cite{P99, P14},
Pedersen~\cite{P00},
Dai~\cite{D01},
Ob{\l}{\'o}j~\cite{O07},
Cox, Hobson and Ob{\l}{\'o}j~\cite{CHO08},
Guo and Zervos~\cite{GZ10},
Glover, Hulley and Peskir~\cite{GHP13},
Galichon, Henry-Labord\`{e}re and Touzi~\cite{GHT14},
Kyprianou and Ott~\cite{KO14},
Gapeev and Rodosthenous~\cite{GR14, GR16},
Rodosthenous and Zervos~\cite{RZ17},
\c{C}a\u{g}lar, Kyprianou and Vardar-Acar~\cite{CKV22},
and D\'{e}camps, Gensbittel and Mariotti~\cite{DGM24},
among others, listed here in rough chronological
order (see also the references therein).
These papers have been motivated by several
applications in mathematical finance and economics 
as well as in quickest detection.

To fix the setting, consider the diffusion $X$ given by
\ben
\di X_t = \mu (X_t) \, \di t + \sigma (X_t) \, \di W_t ,
\quad X_0 = x \in \bbr , \label{Xgen}
\een
where $\mu , \sigma : \bbr \rightarrow \bbr$ are
continuous functions such that $\sigma > 0$.
Let $S$ denote the running maximum process of
$X$, defined by
\ben
S_t = \max \Big\{ s, \max _{0 \leq u \leq t} X_u \Big\} ,
\quad \text{for some } s \geq x .\label{S}
\een
Furthermore, consider the optimal stopping problem
with value function
\ben
v(x,s) = \sup _{\tau \in \tcal}
\EXP _{x,s} \biggl[ \int_0^\tau e^{-rt} f(X_t, S_t) \, \di t
+ \int_0^\tau e^{-rt} g(S_t) \, \di S_t \biggr] ,
\label{v-Lag}
\een
where $\tcal$ is the family of all stopping times,
the discount rate $r \geq 0$ is a constant,
and $f$ and $g > 0$ are continuous functions
satisfying suitable conditions.

By classical results of J.\,L.\,Snell and E.\,B.\,Dynkin, the value function $v$ is the minimal measurable function
$w \geq 0$ such that the process
\be
e^{-rt} w(X_t, S_t) + \int_0^t e^{-ru} f(X_u, S_u) \, \di u
+ \int_0^t e^{-ru} g(S_u) \, \di S_u , \quad  t \geq 0 ,
\ee
is a supermartingale.
The minimal superharmonic characterisation applies to
the optimal stopping of rather general Markov processes
and not merely to the model discussed in the previous
paragraph.
For a wide class of diffusion models, this characterisation
is equivalent to requiring that the value function
satisfies an associated variational inequality in a
suitable sense.
For the purposes of this paper, suppose that
the variational inequality
\ben
\max \bigl\{ \lcal w(x,s) + f(x,s) , \ - w(x,s) \bigr\} = 0 ,
\quad x < s , \label{HJBL}
\een
with the Neumann boundary condition
\ben
w_s (s,s) = - g(s) , \quad s \in \bbr , \label{HJB-BCL}
\een
where
\ben
\lcal w(x,s) = \frac{1}{2} \sigma^2 (x) w_{xx}
(x,s) + \mu (x) w_x (x,s) - r w(x,s) , \label{Lgen}
\een
is solvable in the classical sense.
In this context, the value function $v$ coincides
with the minimal solution to (\ref{HJBL}) and
(\ref{HJB-BCL}).
One way to single out this minimal solution is by
imposing a boundary condition, such as the requirement
that
\ben
\lim _{T \uparrow \infty} e^{-rT} \EXP _{x,s} \bigl[
w(X_T, S_T) \bigr] = 0 , \label{TV-cond}
\een
which is commonly referred to as the transversality
condition.
In many specific problems, however, the minimal solution
to (\ref{HJBL}) and (\ref{HJB-BCL}) can be determined,
without explicitly imposing this condition, by exploiting
structural properties of the problem. 
In particular, when the maximality principle
holds, no additional boundary or transversality condition
is needed. 

Peskir~\cite{P98} solved the optimal stopping problem
defined by (\ref{Xgen})--(\ref{v-Lag}) in the case
$r=0$, $f(x,s) = - c(x)$ and $g(s) = 1$, for a
continuous function $c>0$.
Using the so-called ``principle of smooth fit'', he
derived a continuum of classical solutions to the
variational inequality (\ref{HJBL}) with boundary
condition (\ref{HJB-BCL}), parametrised by
the solutions to a suitable first-order nonlinear ODE
for the optimal stopping problem's free boundary.
Peskir then proved that the solution to the ODE that
corresponds to the minimal solution to (\ref{HJBL})
and (\ref{HJB-BCL}) can be characterised by
the maximality principle, which he stated as:
\ben
\parbox{14cm}{\centering \em
``The optimal stopping boundary $s \mapsto g_* (s)$
for the problem (2.4) \\
is the maximal solution of the differential equation
(3.21) \\
which stays strictly below the diagonal in $\bbr^2$''}
\label{MaxP1}
\een
(see \cite[Section~3.8]{P98}).
In effect, (\ref{MaxP1}) identifies a separatrix
of the ODE \cite[(3.21)]{P98}. 
Peskir's subtle analysis relied on approximating
this separatrix by means of what he termed “bad-good”
solutions, namely, solutions to the ODE \cite[(3.21)]{P98}
that are ``bad'' in the sense that they hit the diagonal of
$\bbr_+^2$ and, at the same time, ``good'' in the sense
that they are not too large. 
A similar phenomenon arises in all other problems
studied in the literature for which the maximality
principle holds.

Later, Peskir~\cite{P14} used the maximality principle to solve
a three-dimensional optimal stopping problem involving
the running minimum process $I$ of $X$, as well as the
running maximum process $S$, and having value function
\be
v(x,i,s) = \sup _{\tau \in \tcal}
\EXP _{x,i,s} \biggl[ S_\tau - I_\tau - \int_0^\tau
f(X_t, I_t, S_t) \, \di t \biggr] .
\ee
Also, Glover, Hulley and Peskir~\cite{GHP13} solved an optimal stopping problem
in which $X$ takes values in $]0, \infty[$, involves
the running minimum process $I$ of $X$, rather than
the running maximum $S$, and has value function
given by
\be
v(x,i) = \inf _{\tau \in \tcal}
\EXP _{x,i} \biggl[ \int_0^\tau \bigl( 1 - 2 L(X_t) L^{-1}(I_t)
\bigr) \, \di t \biggr] ,
\ee
where $L<0$ is a given strictly increasing continuous
function.
In this problem, the observation that
\be
1 - 2 L(x) L^{-1}(i) \leq 0 \quad \Leftrightarrow \quad
x \leq L^\inv \bigl( L(i) / 2 \bigr) =: h(i) ,
\ee
where $L^\inv$ is the inverse function of $L$,
implies that $\bigl\{ (x,i) \in \bbr^2 \mid \
0 < i \leq x < h(i) \bigr\}$ is a subset of the waiting
region.
Furthermore, the optimal stopping boundary is the minimal
solution to the ODE characterising it that stays
strictly above the curve $\bigl( i , h(i) \bigr)$, which
was called ``zero-level curve'' in~\cite{GHP13}.
Effectively, the role of the diagonal in (\ref{MaxP1})
is replaced by the curve defined by the function $h$.
Indeed, one of ``the key novel ingredient[s] revealed in
the solution'' of the problem studied by \cite{GHP13}
``is the replacement of the diagonal and its role in
the maximality principle by a nonlinear curve in the
two-dimensional state space''.

We are thus faced with the following refinement of
the maximality principle.
Consider the optimal stopping problem in
Lagrange form defined by (\ref{Xgen})--(\ref{v-Lag}),
let
\be
\cscr _0 = \bigl\{ (x,s) \mid \ x=s \bigr\} \cup
\bigl\{ (x,s) \mid \ f(x,s) > 0 \bigr\} ,
\ee
and note that $\cscr _0$ is a subset of the waiting
region.
Suppose that the structure of the problem is such that
the principle of smooth fit holds and yields an ODE,
say ($\star$), for (part of) the problem's optimal
boundary.
If the problem's structure does not suggest a
natural boundary condition or asymptotic behaviour
for ($\star$), then
\ben
\parbox{14cm}{\centering \em
the (part of the) optimal stopping boundary characterised
by ($\star$) is the extremal solution to ($\star$) that
remains inside the complement of $\cscr _0$.}
\label{MaxP2}
\een
Here, we use the term extremal, rather than maximal,
to avoid the potential ambiguity arising from the way
the graph of the optimal boundary is plotted.

The expected reward in (\ref{v-Lag}) is formulated in
Lagrange form.
In the present setting, the value function of
optimal stopping problems that are expressed in
Mayer form is defined by
\ben
v (x,s) = \sup _{\tau \in \tcal} \EXP _{x,s} \Bigl[ e^{-r \tau}
R (X_\tau, S_\tau) {\bf 1} _{\{ \tau < \infty \}} \Bigr] ,
\label{v-May}
\een
for a given reward function $R$.
If $R$ is $C^{2,1}$, then an application of
It\^{o}'s formula shows that the optimal stopping
problem defined by (\ref{Xgen}), (\ref{S}) and
(\ref{v-May}) is equivalent to the one given by
(\ref{Xgen})--(\ref{v-Lag}) with
\be
f(x,s) = \lcal R (x,s) \quad \text{and} \quad
g(s) = R_s (s,s) ,
\ee
provided that the associated stochastic integral has
zero expectation.

In this setting, the family of $r$-superharmonic functions
associated with the Markov process $(X,S)$ consists
of all measurable functions $w$ such that the process
$\bigl( e^{-rt} w(X_t, S_t) \bigr)$ is a c\`{a}dl\`{a}g
supermartingale.
Furthermore, the minimal superharmonic
characterisation implies that the value function $v$
is the minimal $r$-superharmonic majorant of $R$,
that is, $v$ is the smallest $r$-superharmonic
function $w$ such that $w \geq R$.
Now, suppose that the variational inequality
\ben
\max \bigl\{ \lcal w(x,s) , \ R(x,s) - w(x,s) \bigr\} = 0 ,
\quad x < s , \label{HJBM}
\een
with the Neumann boundary condition
\ben
w_s (s,s) = 0 , \quad s \in \bbr , \label{HJB-BCM}
\een
where $\lcal$ is given by (\ref{Lgen}), is solvable in
the classical sense.
Then, the value function $v$
coincides with the minimal solution to (\ref{HJBM})
and (\ref{HJB-BCM}).
As in problems formulated in Lagrange form,
this minimal solution can be selected by
using the transversality condition (\ref{TV-cond}).

The maximality principle in the form stated in
(\ref{MaxP2}) has been established in several settings
in which the optimal stopping problem is formulated
directly in Mayer form. 
Representative examples include the problems studied in
\cite{DGM24},
\cite{GZ10} and
\cite{RZ17},
in which $X$ is a geometric Brownian motion, as
well as \cite{P00}, in which $X$ is a general
one-dimensional diffusion.

The purpose of this paper is to present an example
in Mayer form demonstrating that the maximality principle 
does not necessarily hold. 
This finding underscores the need to
rigorously verify that the boundary selected by the
maximality principle indeed coincides with the true
optimal stopping boundary.
In this respect, the problem we solve places the maximality
principle on the same footing as the principle of smooth
fit: it is a highly valuable heuristic tool, but one that
does not always hold.

The example studied in this paper involves the
geometric Brownian motion given by 
\ben
\di X_t = \mu X_t \, \di t + \sigma X_t \, \di W_t ,
\quad X_0 = x > 0 , \label{X}
\een
for some constants $\mu \in \bbr$ and $\sigma \neq 0$,
as well as its running maximum process $S$ that is
defined as in (\ref{S}).
The value function of the optimal stopping problem
is defined in Mayer form by (\ref{v-May}) for $0 < x \leq s$,
with $r>0$ a fixed constant and
\ben
R(x,s) = x^{-1} F(s) - 1 , \quad \text{where}
\quad F(s) = 1 - e^{-s} . \label{F}
\een
We solve this problem by identifying the value function
$v$ with the solution $w$ to the variational inequality
\ben
\max \bigl\{ \lscr w(x,s) , \ x^{-1} F(s) - 1 - w(x,s)
\bigr\} = 0 , \label{HJB}
\een
where
\ben
\lscr w(x,s) = \frac{1}{2} \sigma^2 x^2 w_{xx}
(x,s) + \mu x w_x (x,s) - r w(x,s) , \label{L-op}
\een
subject to the Neumann boundary condition
\ben
w_s (s,s) = 0 , \quad s > 0 , \label{HJB-BC}
\een
and the transversality condition (\ref{TV-cond}). 
In deriving $w$, we use the principle of smooth fit. 

The value function of the equivalent optimal stopping
problem in Lagrange form is
\be
\sup _{\tau \in \tcal} \EXP _{x,s} \biggl[ \int_0^{\tau}
e^{-rt} \lscr R(X_t,S_t) \, \di t + \int_0^\tau e^{-rt} R_s
(S_t, S_t) \, \di S_t \biggr] .
\ee
The reward function $R$ is such that
\be
\lscr R(x,s) \begin{cases} < 0 , & \text{for all }
x \in \mbox{} \bigl] 0, G(s) \bigr[ , \\ > 0 , &
\text{for all } x \in \mbox{} \bigl] G(s) , s \bigr[
, \end{cases} \quad \text{and} \quad
R_s (s,s) > 0 \text{ for all } s > 0 ,
\ee
where $G : \bbr_+ \rightarrow \bbr_+$ is
the strictly increasing function defined by (\ref{G}).
Consequently, the domain
\be
\cscr _0 = \big\{ (x,s) \in \bbr_+^2 \mid \ G(s) < x
\leq s \big\} ,
\ee
defined by the ``zero-level curve'' $x = G(s)$,
must be a subset of the optimal stopping problem's
waiting region $\wscr$, otherwise, the variational
inequality (\ref{HJB}) would not hold.
In light of this observation and the structure of
the optimal stopping problem, we seek
a strictly increasing free-boundary function
$H$ that separates the waiting region $\wscr$ from
the stopping region $\sscr$ and satisfies
\ben
H(0) = 0 \quad \text{and} \quad 0 < H(s) < G(s)
\text{ for all } s > 0 . \label{H-reqs0}
\een
The boundary condition $H(0) = 0$ is not a
consequence of the heuristic arguments above.
Instead, it is derived from the subsequent rigorous
analysis.
Specifically, the stopping and waiting regions are
given by
\be
\sscr = \Big\{ (x,s) \in \bbr_+^2 \mid \ 0 < x \leq
H(s) \Big\} \quad \text{and} \quad
\wscr = \Big\{ (x,s) \in \bbr_+^2 \mid \ H(s) < x \leq
s \Big\} .
\ee

\begin{figure}[!tbp]
\centering
\includegraphics[width=150mm]{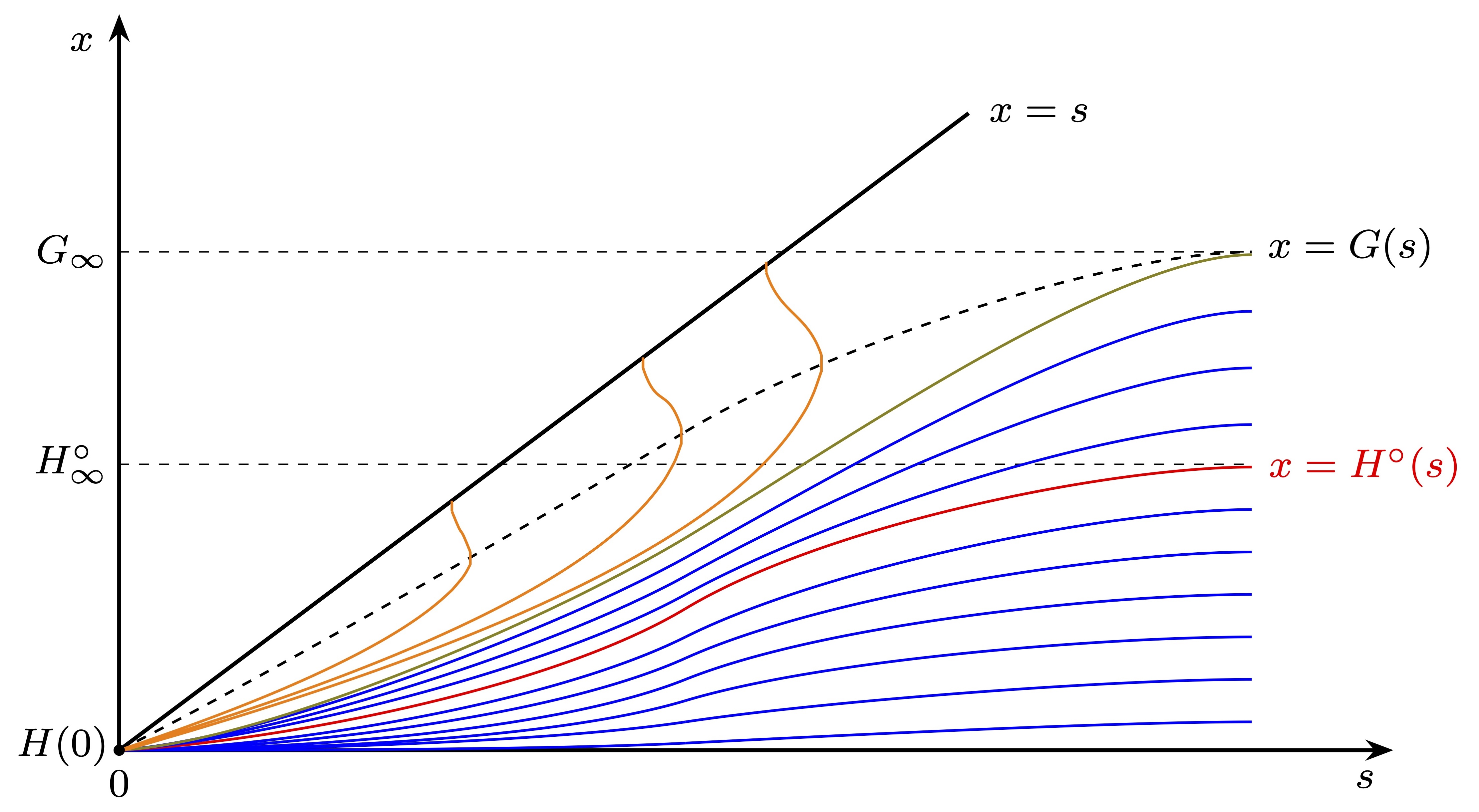}
\caption{\footnotesize Illustration of possible solutions
$H$ to the ODE (\ref{H-ODE}).
The red curve corresponds to the solution identifying
the optimal stopping boundary. 
The green curve represents the solution that would arise
in the context of the maximality principle as stated in
(\ref{MaxP2}), with the dotted black line indicating the
``zero-level curve'' $x = G(s)$.
The blue curves represent other solutions to (\ref{H-ODE})
that satisfy (\ref{H-reqs0}).}
\label{fig-H}
\end{figure}

We show that the boundary condition
(\ref{HJB-BC}) and an application of the principle
of smooth fit yield the first-order ODE
\ben
\dot{H} (s) = \frac{\dot{F} (s) \Bigl( (n+1)
\bigl( H(s) / s \bigr)^{n-m} - (m+1) \Bigr) H(s)}
{-mn \bigl( G(s) - H(s) \bigr) \Bigl( 1 - \bigl(
H (s) / s \bigr) ^{n-m} \Bigr)} \label{H-ODE}
\een
for the free-boundary function $H$.
Figure~\ref{fig-H} illustrates the continuum of solutions
to this ODE.
In Theorems~\ref{thm:w-VI} and~\ref{thm:main},
we prove that the optimal stopping boundary
$s \mapsto H^\circ (s)$ is the unique solution to
the ODE (\ref{H-ODE}) that is associated with a
solution $w$ to the variational inequality (\ref{HJB})
and the boundary condition (\ref{HJB-BC})
that satisfies the transversality condition
(\ref{TV-cond}).
In particular, we prove that
(a) there exists a continuum of solutions $H$
to the ODE (\ref{H-ODE}) satisfying (\ref{H-reqs0})
and such that $H^\circ < H < G$, while
(b) each of these solutions corresponds to a
candidate function $w$ that fails to satisfy the
variational inequality (\ref{HJB}) because there
are points $0<x<s$ such that $w(x,s) < 0$.
This analysis shows that the optimal stopping boundary
lies strictly below the solution to the associated
ODE that would be selected by the maximality principle
(see red and green curves in Figure~\ref{fig-H},
respectively).

\section{The optimal stopping problem}
\label{sec:solution}

We fix a filtered probability space $\bigl( \Omega, \fcal,
(\fcal_t), \PR \bigr)$ satisfying the usual conditions
and supporting a standard one-dimensional
$(\fcal_t)$-Brownian motion $W$.
We denote by $\tcal$ the set of all $(\fcal_t)$-stopping
times.\footnote{
In this paper, we adopt a strong rather than a weak
formulation.
In the introduction, we used the notation $\EXP_{x,s}$
to emphasise the dependence of the process $(X,S)$
on its initial condition $(X_0, S_0) = (x,s)$, rather than
to indicate that expectations are computed with respect
to a probability measure $\PR_{x,s}$.
In the remainder of the paper, we omit the subscript $(x,s)$
from $\EXP_{x,s}$ for notational simplicity.
}

The solution to the optimal stopping problem defined
by (\ref{v-May}) with $S$, $X$ and $R$ given by
(\ref{S}), (\ref{X}) and (\ref{F}) involves the general
solution to the ODE
\ben
\frac{1}{2} \sigma^2 x^2 f''(x) + \mu x f'(x) - r f(x) = 0
, \label{ODE}
\een
which is given by
\ben
f(x) = A x^n + B x^m , \label{ODE-sol}
\een
for some $A, B \in \bbr$.
Here, the constants $m < 0 < n$ are the solutions to the
quadratic equation
\ben
\frac{1}{2} \sigma^2 k^2 + \biggl( \mu - \frac{1}{2} \sigma^2
\biggr) k - r = 0 , \label{mn-eq}
\een
which are given by
\ben
m, n = - \frac{\mu - \frac{1}{2} \sigma^2}{\sigma^2} \mp
\sqrt{\left( \frac{\mu - \frac{1}{2} \sigma^2}{\sigma^2}
\right)^2 + \frac{2 r}{\sigma^2}} . \label{mn}
\een

\begin{ass} \label{A}  {\rm
The constants $r > 0$, $\mu \in \bbr$ and $\sigma \neq 0$
are such that
\ben
m+1 < 0 \ \Leftrightarrow \ \sigma^2 < r + \mu
\quad \text{and} \quad
m+n+1 \geq 0 \ \Leftrightarrow \ \sigma^2 \geq \mu .
\label{A:eqn}
\een
} \end{ass}

The first of the conditions in (\ref{A:eqn}) guarantees
that the value function $v$ is finite and that an optimal
stopping time indeed exists.
In particular, it implies that
\ben
\lim _{\ell \uparrow \infty} \EXP \Bigl[
e^{-r \tau_\ell} \bigl( X_{\tau_\ell}^{-1}
F(S_{\tau_\ell}) - 1 \bigr) \Bigr] \leq
\lim _{\ell \uparrow \infty} \EXP \bigl[
e^{-r \tau_\ell} X_{\tau_\ell}^{-1} \bigr]
= 0 , \label{TV-F}
\een
for any sequence $(\tau_\ell)$ of bounded
$(\fcal _t)$-stopping times such that $\lim
_{\ell \uparrow \infty} \tau_\ell = \infty$.
We assume that the second one also holds
true because it implies that
\ben
0 < \frac{(m+1) (n+1)}{mn} \leq 1, \label{frac<1}
\een
which will simplify our analysis.

We will solve the optimal stopping problem we consider
by identifying the value function
$v$ with the solution $w$ to the variational inequality
(\ref{HJB}) with the Neumann boundary condition
(\ref{HJB-BC}) that satisfies the transversality condition
(\ref{TV-cond}).
A solution $w$ to the variational inequality (\ref{HJB}) 
partitions the problem's state space into
the sets
\begin{align}
\sscr & = \Big\{ (x,s) \in \bbr_+^2 \mid \ 0 < x \leq
s \text{ and } w(x,s) = x^{-1} F(s) - 1
\Big\} \label{Sregion} \\
\text{and} \quad
\wscr & = \Big\{ (x,s) \in \bbr_+^2 \mid \ 0 < x \leq
s \text{ and } w(x,s) > x^{-1} F(s) - 1
\Big\} ,  \label{Wregion}
\end{align}
which are candidates for the so-called stopping
and waiting regions, respectively.

In view of the definitions (\ref{mn}) of $m$ and $n$,
we can see that
\begin{align}
\frac{1}{2} \sigma^2 x^2
\frac{\partial^2 \bigl( x^{-1} F(s) - 1 \bigr)}{\partial x^2}
+ \mu x & \frac{\partial \bigl(  x^{-1}F(s) - 1 \bigr)}
{\partial x} - r \bigl( x^{-1} F(s) - 1 \bigr) > 0
\nonumber \\
& \Leftrightarrow \quad
(\sigma^2 - \mu - r) x^{-1} F(s) + r > 0 \nonumber \\
& \Leftrightarrow \quad
x > \frac{(m+1) (n+1)}{mn} F(s) =: G(s) . \label{G}
\end{align}
These equivalences imply that any classical solution $w$
to (\ref{HJB}) is such that
\be
\big\{ (x,s) \in \bbr_+^2 \mid \ G(s) < x \leq s \big\}
\subseteq \wscr
\quad \text{and} \quad
\sscr \subseteq \big\{ (x,s) \in \bbr_+^2 \mid \ 0 <
x \leq G(s) < s \big\} .
\ee
Motivated by this observation, we will look for a solution
to (\ref{HJB}) that is associated with the stopping and
waiting regions
\ben
\sscr = \bigl\{ (x,s) \in \bbr_+^2 \mid \ 0 < x \leq H(s) \bigr\}
\quad \text{and} \quad
\wscr = \bigl\{ (x,s) \in \bbr_+^2 \mid \ H(s) < x \leq s \bigr\}
, \label{SW}
\een
for some strictly increasing free-boundary function
$H$ such that
\ben
H(0) = 0 \quad \text{and} \quad 0 < H(s) < G(s)
\text{ for all } s > 0 . \label{H-reqs}
\een
Accordingly, we will look for a solution to (\ref{HJB})
that is of the form
\ben
w(x,s) = \begin{cases} x^{-1} F(s) - 1 , & \text{if }
x \in \mbox{} ]0, H(s)] , \\  A(s) x^n + B (s) x^m ,
& \text{if } x \in \mbox{} ]H (s), s] , \end{cases}
\label{w}
\een
for some functions $A$ and $B$, because $w(\cdot , s)$
should satisfy the ODE (\ref{ODE}) in the interior of
the waiting region $\wscr$.

For future reference, we note that the function
$G: \mbox{} ]0, \infty[ \mbox{} \to \bbr$ defined
by (\ref{G}) is strictly increasing and strictly concave,
\begin{gather}
G(s) < s \quad \text{for all } s > 0 , \label{G(s)<s} \\
G_0 :=  \lim _{s \downarrow 0} G(s) = 0
\quad \text{and} \quad
G_\infty := \lim _{s \uparrow \infty} G(s) =
\frac{(m+1) (n+1)}{mn} \in \mbox{} ] 0,1] ,
\label{Ginf}
\end{gather}
by virtue of the inequalities (\ref{frac<1}).

\section{The free-boundary function}
\label{sec:H}

To determine the functions $A$ and $B$, as well as the
free-boundary function $H$, appearing in (\ref{w}), we
first note that the boundary condition (\ref{HJB-BC})
gives rise to the equation
\ben
\dot{A}(s) s^n + \dot{B}(s) s^m = 0 , \quad \text{for }
s > 0 . \label{W2BC}
\een
In view of the regularity of the optimal stopping problem's
reward function, we expect that the so-called ``principle
of smooth fit'' should hold true.
Accordingly, we require that $w(\cdot, s)$ should be $C^1$
along the free-boundary point $H(s)$, for all $s>0$.
This requirement yields the system of equations
\begin{gather}
A(s) H^n (s) + B(s) H^m (s) = H^{-1} (s) F(s) - 1
\nonumber \\
\text{and} \quad
n A(s) H^n (s) + m B(s) H^m (s) = - H^{-1} (s) F(s) ,
\nonumber
\end{gather}
which is equivalent to
\ben
A(s) = \frac{m H(s) - (m+1) F(s)}{(n-m) H^{n+1} (s)}
\quad \text{and} \quad
B(s) = \frac{(n+1) F(s) - n H(s)}{(n-m) H^{m+1} (s)} .
\label{AnB}
\een
Differentiating these expressions with respect to $s$
and substituting the results for $\dot{A}$ and $\dot{B}$
in (\ref{W2BC}), we can see that $H$ should satisfy
the ODE (\ref{H-ODE}), where $G$ is defined by
(\ref{G}).

A function $H$ satisfies the ODE (\ref{H-ODE})
in the domain
\be
\dscr_H = \bigl\{ (s,h) \in \bbr_+^2 \mid \ 0 < h < s
\text{ and } h \neq G(s) \bigr\}
\ee
if and only if the function $Q$ defined by
\be
Q(s) = \frac{H(s)}{F (s)} , \quad \text{for } s > 0 ,
\ee
satisfies the ODE
\ben
\dot{Q} (s) = \frac{\dot{F}(s)}{F (s)} Q(s)
\qcal \bigl( s, Q(s) \bigr)
\quad \Leftrightarrow \quad
\frac{\di \ln Q(s)}{\di s} = \qcal \bigl( s, Q(s) \bigr)
\frac{\di \ln F(s)}{\di s} \label{Q-ODE}
\een
in the domain
\be
\dscr_Q = \biggl\{ (s,q) \in \bbr_+^2 \mid \ 0 <
q < \frac{1}{\gamma (s)} \text{ and } q \neq
G_\infty \biggr\} ,
\ee
which corresponds to $\dscr_H$, where
\ben
\qcal (s,q) = \frac{\Bigl( \frac{n+1}{n} - q \Bigr)
\bigl( \gamma (s) q \bigr)^{n-m} + q - \frac{m+1}{m}} 
{\bigl( G_\infty - q \bigr) \Bigl( 1 - \bigl(
\gamma (s) q \bigr)^{n-m} \Bigr)}
\quad \text{and} \quad
\gamma (s) = \frac{F(s)}{s} .
\label{qcal-gam}
\een
For future reference, we note
\ben
\lim _{s \downarrow 0} \gamma (s) = 1 , \quad
\dot{\gamma} (s) < 0 \text{ for all } s>0
\quad \text{and} \quad
\lim _{s \uparrow \infty} \gamma (s) = 0 .
\label{gamma-props}
\een

\begin{lem} \label{lem:zeta}
There exists a point $q_\dagger \in \mbox{}
\bigl] 0, \frac{m+1}{m} \bigr[$ and a continuous function
$\zeta : \bbr_+ \rightarrow \bbr_+$ such that
\begin{gather}
\lim _{s \downarrow 0} \zeta (s) = q_\dagger ,
\quad \dot{\zeta} (s) > 0 \text{ for all } s>0 ,
\quad \lim _{s \uparrow \infty} \zeta (s) =
\frac{m+1}{m} \nonumber \\
\text{and} \quad
\Bigl( \frac{n+1}{n} - q \Bigr) \bigl( \gamma (s)
q \bigr)^{n-m} + q - \frac{m+1}{m}
\begin{cases}
< 0 & \text{for all } s>0 \text{ and } q < \zeta (s) , \\
> 0 & \text{for all } s>0 \text{ and } q > \zeta (s) .
\end{cases} \nonumber
\end{gather}
\end{lem}
\noindent {\bf Proof.}
We first note that
\ben
\biggl( \frac{n+1}{n} - q \biggr) \bigl( \gamma
(s) q \bigr)^{n-m} + q - \frac{m+1}{m}
\geq \frac{n-m}{-mn} > 0 \label{qnum>0}
\een
for all $(q,s)$ such that $s > 0$ and $\frac{n+1}{n}
\leq q \leq \frac{1}{\gamma (s)}$.
We next define
\be
\theta (q) = q \biggl(
\frac{\frac{n+1}{n} -q}{\frac{m+1}{m} - q}
\biggr) ^{1/(n-m)} , \quad \text{for }
q \in \mbox{} \bigl] 0, (m+1)/m \bigr[ ,
\ee
and we observe that
\be
\lim _{q \downarrow 0} \theta (q) = 0
, \quad \dot{\theta} (q) > 0 \text{ for all } q \in \bigl] 0, (m+1)/m \bigr[ 
\quad \text{and} \quad \lim _{q \uparrow \frac{m+1}{m}}
\theta (q) = \infty .
\ee
If we define $\zeta (s) = \theta^\inv \bigl( 1 /
\gamma (s) \bigr)$, where $\theta^\inv$ is the
inverse function of $\theta$ and $\gamma$ is
defined by (\ref{qcal-gam}), then all of the claims
of the lemma follow from (\ref{gamma-props})
and (\ref{qnum>0}) for $q_\dagger =
\theta^\inv (1)$.
\mbox{}\hfill$\Box$
\bigskip

In view of the definition (\ref{qcal-gam}) of $\qcal$
and the previous result, we can see that
\begin{gather}
\qcal (s,q) > 0 \quad \text{for all } (s,q) \in \dscr
_Q^+ , \label{qcal-props1} \\
\qcal (s,q) < 0 \quad \text{for all } (s,q) \in \dscr
_Q^{l-} \cup \dscr _Q^{u-} ,
\label{qcal-props2} \\
\qcal \bigl( s, \zeta (s) \bigr) = 0 \quad
\text{for all } s > 0 , \label{qcal-props3} \\
\lim _{q \uparrow G_\infty} \qcal (s,q) = \infty
\quad \text{and} \quad
\lim _{q \downarrow G_\infty} \qcal (s,q)
= \lim _{q \uparrow \frac{1}{\gamma (s)}}
\qcal (s,q) = -\infty \quad \text{for all }
s > 0 , \label{qcal-props4}
\end{gather}
where
\begin{gather}
\dscr _Q^+ = \Bigl\{ (s,q) \in \bbr_+^2 \mid
\ \zeta (s) < q < G_\infty \Bigr\} , \label{DQ+} \\
\dscr _Q^{l-} = \Bigl\{ (s,q) \in \bbr_+^2 \mid \ 0
< q < \zeta (s) \Bigr\} , \quad
\dscr _Q^{u-} = \biggl\{ (s,q) \in \bbr_+^2 \mid \ 
G_\infty < q < \frac{1}{\gamma (s)} \biggr\} .
\label{DQ-}
\end{gather}
Given any $\varepsilon \in
\mbox{} ]0, 1[$, the restriction of $\qcal$ in the
domain
\be
\biggl\{ (s,q) \in \bbr_+^2 \mid \ \varepsilon < s
< \frac{1}{\varepsilon} \text{ and } \Bigl(
0 < q < G_\infty {-} \varepsilon \text{ or }
G_\infty {+} \varepsilon < q <
\frac{1}{\gamma (s)} {-} \varepsilon \Bigr)
\biggr\}
\ee
is Lipschitz continuous.
Therefore, given any $(s_0 , q_0) \in \dscr_Q$,
there exist
\ben
\underline{s} (s_0 , q_0) \in [0, s_0[
\quad \text{and} \quad
\overline{s} (s_0 , q_0) \in \mbox{}]
s_0, \infty] \label{undover-s0}
\een
such that the ODE (\ref{Q-ODE}) with initial
condition $Q (s_0) = q_0$ has a unique solution
$Q_{(s_0, q_0)}$ satisfying
\begin{gather}
\lim _{s \downarrow \underline{s}} Q_{(s_0, q_0)}
(s) \begin{cases} \in \mbox{} ]0, 1[ , & \text{if }
\underline{s} (s_0, q_0) = 0 , \\ =
\frac{1}{\gamma (\underline{s})} , & \text{if }
\underline{s} (s_0, q_0) > 0 , \end{cases}
\label{undover-s1} \\
\text{and } \quad
\lim _{s \uparrow \overline{s}} Q_{(s_0, q_0)}
(s) \begin{cases} \in \{ 0 , G_\infty \} , & \text{if }
\overline{s} (s_0, q_0) < \infty , \\ \in \mbox{}
]0, \infty[ , & \text{if } \overline{s} (s_0, q_0)
= \infty , \end{cases} \label{undover-s2}
\end{gather}
(see Piccinini, Stampacchia and Vidossich~\cite[Theorems
I.1.4 and~I.1.5]{PSV}).
Furthermore,
\begin{align}
\text{if } q_0^1 < q_0^2 , \quad & \text{then }
Q_{(s_0, q_0^1)} (s) < Q_{(s_0, q_0^2)} (s)
\nonumber \\
& \text{for all } s \in \mbox{} \bigl]
\underline{s} (s_0, q_0^1) {\vee}
\underline{s} (s_0, q_0^2) , \,
\overline{s} (s_0, q_0^1) {\wedge}
\overline{s} (s_0, q_0^2) \bigr[ ,
\label{Q-order}
\end{align}
by virtue of the uniqueness of solutions to
the ODE (\ref{Q-ODE}) in $\dscr_Q$.

\begin{figure}[!tbp]
\centering
\includegraphics[width=120mm]{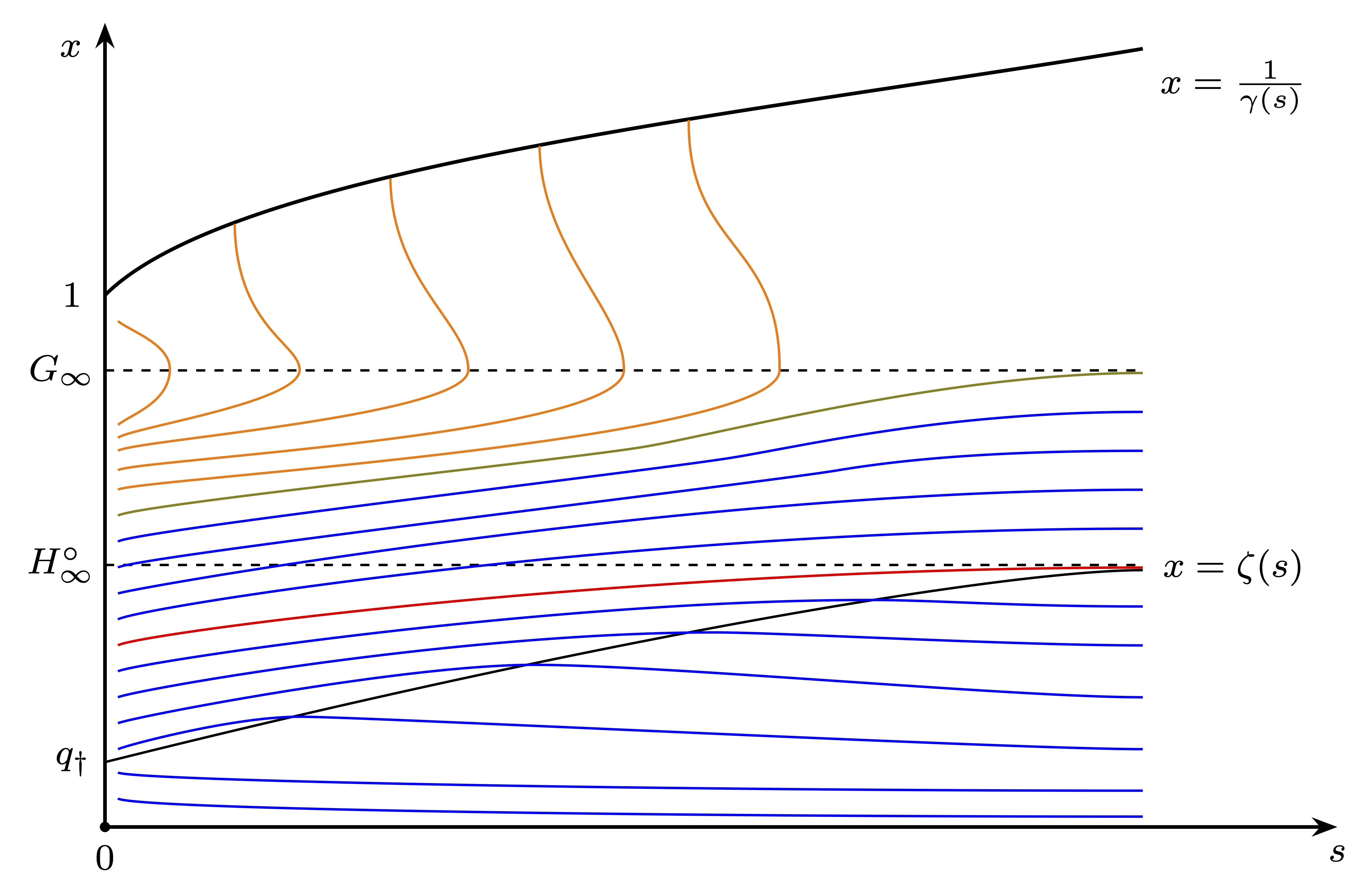}
\caption{\footnotesize Illustration of possible solutions
to the ODE (\ref{Q-ODE}).
The level $G_\infty$ is defined as in (\ref{Ginf}),
while $H_\infty^\circ = \frac{m+1}{m}$.
The red curve represents the solution $Q$ to
(\ref{Q-ODE}) that is associated with the optimal
stopping boundary $H$.
The green curve represents the separatrix that separates
solutions $Q$ to (\ref{Q-ODE}), such as the ones associated
with the orange curves, from solutions $Q$ to
(\ref{Q-ODE}) that correspond to solutions $H$
to the ODE (\ref{H-ODE}) satisfying (\ref{H-reqs0})
(blue, green and red curves).}
\label{fig-Q}
\end{figure}

The following result, which is illustrated by Figure~\ref{fig-Q},
presents a study of the ODE (\ref{Q-ODE}).

\begin{thm} \label{prop:Q}
Suppose that Assumption~\ref{A} holds true and
consider the domains $\dscr _Q^+$, $\dscr _Q^{l-}$,
$\dscr _Q^{u-}$ defined by (\ref{DQ+}) and
(\ref{DQ-}), as well as the points $\underline{s}
= \underline{s} (s_0, q_0) < \overline{s} (s_0, q_0)
= \overline{s}$ associated with each $(s_0, q_0)
\in \dscr_Q$ and (\ref{undover-s0})--(\ref{undover-s2}).
The following statements hold true.
\smallskip

\noindent {\rm (I)}
Given any point $(s_0 , q_0) \in \dscr _Q^{u-}$, the
ODE (\ref{Q-ODE}) with initial condition
$Q (s_0) = q_0$ has a unique solution such that
\be
\lim _{s \downarrow \underline{s}} Q(s)
\begin{cases} \in \mbox{} ]G_\infty, 1[ , & \text{if }
\underline{s} = 0 , \\ = \frac{1}{\gamma (\underline{s})}
, & \text{if } \underline{s} > 0 , \end{cases} \quad
\dot{Q} (s) < 0 \text{ for all } s \in \mbox{}
\bigl] \underline{s} , \overline{s} \bigr[
\quad \text{and} \quad
\lim _{s \uparrow \overline{s}} Q(s) \geq G_\infty .
\ee

\noindent {\rm (II)}
Given any point $(s_0 , q_0) \in \dscr _Q^{l-}$, the
ODE (\ref{Q-ODE}) with initial condition
$Q (s_0) = q_0$ has a unique solution such that
\be
\overline{s} (s_0, q_0) = \infty 
\quad \text{and} \quad
\dot{Q} (s) < 0 \text{ for all } s > s_0 .
\ee

\noindent {\rm (III)}
Given any $q_\infty \in \mbox{} ]0, G_\infty]$,
there exists a unique solution to the ODE
(\ref{Q-ODE}) in $\dscr _Q^{l-} \cup \dscr _Q^+$
such that
\be
\lim _{s \downarrow 0} Q(s) \in \mbox{} ]0,
G_\infty[ \quad \text{and} \quad
\lim _{s \uparrow \infty} Q(s) = q_\infty .
\ee
Furthermore, if $q_\infty \in \bigl[ \frac{m+1}{m}
, G_\infty \bigr]$, then this solution takes
values in $\dscr _Q^+$, in particular,
\be
\lim _{s \downarrow 0} Q(s) \in \bigl[ q_\dagger ,
G_\infty \bigr[ \quad \text{and} \quad
\dot{Q} (s) > 0 \text{ for all } s > 0 ,
\ee
where $q_\dagger \in \mbox{} \bigl] 0,
\frac{m+1}{m} \bigr[$ is as in
Lemma~\ref{lem:zeta}.
\end{thm}
\noindent {\bf Proof.}
The claims in (I) follow immediately from
(\ref{qcal-props2}) and the last two limits in
(\ref{qcal-props4}).

To prove~(II), fix any point $(s_0 , q_0) \in \dscr
_Q^{l-}$.
If $Q$ is the solution to the ODE (\ref{Q-ODE})
with initial condition $Q(s_0) = q_0$, then
(\ref{qcal-props2}) implies that $\dot{Q} (s) < 0$
for all $s \in \mbox{} \bigl] s_0, \overline{s}
(s_0, q_0) \bigr[$.
To show that $\overline{s} (s_0, q_0) = \infty$,
we write $\qcal (s,q) = N(s,q) / D(s,q)$ and we
note that the calculation
\begin{align}
\frac{\partial N(s,q)}{\partial q} & = \biggl(
\frac{(n+1) (n-m)}{n} - (n-m+1) q \biggr)
\frac{1}{q} \bigl( \gamma (s) q \bigr)^{n-m}
+ 1 \nonumber \\
& > 1 \quad \text{for all } s > 0 \text{ and } q < 1
, \nonumber
\end{align}
the fact that $G_\infty \leq 1$, the inequality
$\partial D(s,q) / \partial q < 0$ (which holds true
for all $q \in \mbox{} ]0, G_\infty[$) and the
assumption that $m+1 < 0$ (see (\ref{A:eqn}))
imply that
\be
\qcal (s,q) \geq - \frac{m+1} 
{m \bigl( G_\infty - q_0 \bigr) \Bigl( 1 - \bigl(
\gamma (s_0) q_0 \bigr)^{n-m} \Bigr)}
=: \beta (s_0 , q_0) \quad \text{for all }
s \geq s_0 \text{ and } q \leq q_0 .
\ee
In view of this inequality and the ODE (\ref{Q-ODE}),
we can see that
\be
\int _{s_0}^s \di \ln Q(u) = \int _{s_0}^s \qcal
\bigl( u, Q(u) \bigr) \, \di \ln F(u) \geq \beta
(s_0 , q_0) \int _{s_0}^s \di \ln F(u) .
\ee
It follows that
\be
Q(s) \geq \frac{q_0}{F^{\beta (s_0 , q_0)} (s_0)}
F^{\beta (s_0 , q_0)} (s) > 0 \quad \text{for all }
s > s_0 ,
\ee
which implies that $\overline{s} (s_0, q_0)
= \infty$.

To establish (III), we fix any $q_\infty \in \mbox{}
]0 , G_\infty[$ and any $\varepsilon > 0$ such that
$\bigl] (1 {-} \varepsilon) q_\infty ,
(1 {+} \varepsilon) q_\infty \bigr[
\mbox{} \subseteq \mbox{} ]0, G_\infty[$.
Also, we choose any $\widetilde{s}_\varepsilon
= \widetilde{s} _\varepsilon (q_\infty) > 0$
such that
\be
F^{-\ell (q_\infty) + \varepsilon} (s) \in \mbox{} ]1 {-}
\varepsilon , 1 {+} \varepsilon[
\quad \text{and} \quad
F^{-\ell (q_\infty) - \varepsilon} (s) \in \mbox{} ]1 {-}
\varepsilon , 1 {+} \varepsilon[ \quad
\text{for all } s \geq \widetilde{s}_\varepsilon ,
\ee
where
\be
\ell (q) = \frac{q - \frac{m+1}{m}}
{G_\infty - q} .
\ee
Such a choice implies that
\be
\frac{F^{\ell (q_\infty) - \varepsilon} (s)}
{F^{\ell (q_\infty) - \varepsilon}
(\widetilde{s}_\varepsilon)} \in \mbox{}
]1 {-} \varepsilon , 1 {+} \varepsilon[
\quad \text{and} \quad
\frac{F^{\ell (q_\infty) + \varepsilon} (s)}
{F^{\ell (q_\infty) + \varepsilon}
(\widetilde{s}_\varepsilon)} \in \mbox{}
]1 {-} \varepsilon , 1 {+} \varepsilon[ \quad
\text{for all } s \geq \widetilde{s}_\varepsilon
\ee
because
\be
\frac{F^\alpha (s)}
{F^\alpha (\widetilde{s}_\varepsilon)} \in
\begin{cases}
\bigl] F^{-\alpha} (\widetilde{s}_\varepsilon) ,
1 \bigr] , & \text{if } \alpha < 0 , \\
\bigl[ 1 , F^{-\alpha} (\widetilde{s}_\varepsilon)
\bigr[ , & \text{if } \alpha > 0 ,
\end{cases} \quad \text{for all } s \geq
\widetilde{s}_\varepsilon .
\ee
Next, we note that $\lim _{s \uparrow \infty}
\qcal (s,q_\infty) = \ell (q_\infty)$ (see
(\ref{qcal-gam}), (\ref{gamma-props}) and
the definition of $\ell$ above).
%
In view of this observation, we choose a constant
$\delta = \delta(\varepsilon) \in \mbox{} ]0,\varepsilon]$
and any
$s_\varepsilon = s_\varepsilon (q_\infty)
\geq \widetilde{s} _\varepsilon (q_\infty)$
such that
\ben
\qcal (s, q) \in \mbox{} \bigl] \ell (q_\infty) {-}
\varepsilon , \ell (q_\infty) {+} \varepsilon
\bigr[ \quad \text{for all } s \geq s_\varepsilon
\text{ and } q \in \mbox{} \bigl] (1 {-} \delta)
q_\infty , (1 {+} \delta) q_\infty \bigr[ ,
\label{Qlim-eps}
\een
and $F^{-\ell(q_\infty) + \varepsilon} (s_\varepsilon),
F^{-\ell(q_\infty) - \varepsilon} (s_\varepsilon) \in
\mbox{} ]1 {-} \delta , 1 {+} \delta[$.

If $Q$ is the solution to the ODE (\ref{Q-ODE})
with initial condition $Q(s_\varepsilon) = q_\infty$,
then (\ref{Qlim-eps}) and the observation that
\be
\int _{s_\varepsilon}^s \di \ln Q(u) =
\int _{s_\varepsilon}^s \qcal \bigl( u, Q(u) \bigr)
\, \di \ln F(u) \begin{cases}
\geq \bigl( \ell (q_\infty) - \varepsilon \bigr) \int
_{s_\varepsilon}^s \di \ln F(u) , \\
\leq \bigl( \ell (q_\infty) + \varepsilon \bigr) \int
_{s_\varepsilon}^s \di \ln F(u) ,
\end{cases}
\ee
imply that
\be
(1 {-} \delta) q_\infty <
q _\infty\frac{F^{\ell (q_\infty) - \varepsilon} (s)}
{F^{\ell (q_\infty) - \varepsilon} (s_\varepsilon)}
\leq Q(s) \leq
q_\infty \frac{F^{\ell (q_\infty) + \varepsilon} (s)}
{F^{\ell (q_\infty) + \varepsilon} (s_\varepsilon)}
< (1 {+} \delta) q_\infty .
\ee
Therefore,
\be
\lim _{s \uparrow \infty} Q(s) \in \mbox{} \bigl] (1 {-} \delta) q_\infty ,  (1 {+} \delta) q_\infty \bigr[ \, \subseteq \mbox{}
\bigl] (1 {-} \varepsilon) q_\infty ,  (1 {+}
\varepsilon) q_\infty \bigr[ .
\ee
Furthermore, (\ref{qcal-props1}), (\ref{qcal-props2})
and (\ref{qcal-props4}) imply that this solution is
well defined for all $s>0$ and such that
$\lim _{s \downarrow 0} Q(s) \in \mbox{} ]0,
G_\infty[$.

Fix any $\overline{s} > 0$.
The analysis above establishes that, given any
$q_\infty \in \mbox{} ]0, G_\infty[$ and any
$\varepsilon > 0$  such that $\bigl] (1 {-} \varepsilon)
q_\infty , (1 {+} \varepsilon) q_\infty \bigr[ \mbox{}
\subseteq \mbox{} ]0, G_\infty[$, there exists a
point $\overline{q} = \overline{q} (\varepsilon,
q_\infty)$ such that the ODE (\ref{Q-ODE}) with
initial condition $Q(\overline{s}) = \overline{q}$
has a unique solution $Q_{(\overline{s}, \overline{q})}$
such that
\be
\lim _{s \downarrow 0} Q_{(\overline{s}, \overline{q})}
(s) \in \mbox{} ]0, G_\infty[ \quad \text{and} \quad
\lim _{s \uparrow \infty} Q_{(\overline{s}, \overline{q})}
(s) \in \mbox{} \bigl] (1 {-} \varepsilon) q_\infty
,  (1 {+} \varepsilon) q_\infty \bigr[ .
\ee
This observation, the fact that $\varepsilon > 0$
can be arbitrarily small and the continuous dependence
of the solution to an ODE with respect to its
initial conditions imply all of the claims in
part~(III) for $q_\infty \in \mbox{} ]0 , G_\infty[$.

To proceed further, we parametrise the solutions
derived in the previous paragraph by their limiting
value $q_\infty \in \mbox{} ]0 , G_\infty[$ and we
write $Q (\cdot; q_\infty)$ instead of $Q$.
Furthermore, we define
\be
Q^\circ (s) = \lim _{q_\infty \uparrow G_\infty}
Q(s; q_\infty) < G_\infty , \quad
\text{for } s \geq 0 .
\ee
The strict inequality here is an immediate
consequence of the first limit in (\ref{qcal-props4}).
In view of (\ref{Q-order}), we use the monotone
or the dominated convergence theorems to obtain
\begin{align}
Q^\circ (s_2) & = Q^\circ (s_1) + \lim
_{q_\infty \uparrow G_\infty} \int _{s_1}^{s_2}
\frac{\dot{F}(u)}{F (u)} Q(u; q_\infty)
\qcal \bigl( u, Q(u; q_\infty) \bigr) \, \di u
\nonumber \\
& = Q^\circ (s_1) + \int _{s_1}^{s_2}
\frac{\dot{F}(u)}{F (u)} Q^\circ (u) \qcal \bigl(
u, Q^\circ (u) \bigr) \, \di u \quad
\text{for all } s_1 < s_2 . \nonumber
\end{align}
It follows that $Q^\circ$ is a solution to the
ODE (\ref{Q-ODE}) such that $\lim
_{s \uparrow \infty} Q^\circ (s) = G_\infty$.
\mbox{}\hfill$\Box$
\bigskip

In the following result, we consider
only solutions to the ODE (\ref{H-ODE}) that can
be identified with the optimal stopping problem's
free-boundary function $H$, namely, solutions
that satisfy (\ref{H-reqs}) (see also
Figure~\ref{fig-H}).

\begin{cor} \label{cor:H-ODE}
Suppose that Assumption~\ref{A} holds true.
Given any $H_\infty \in \mbox{} ]0 , G_\infty]$,
the ODE (\ref{H-ODE}) has a unique solution
$H$ such that
\begin{gather}
0 < H(s) < G(s) \quad \text{and} \quad
\dot{H} (s) > 0 \quad \text{for all } s > 0 ,
\nonumber \\
\lim _{s \downarrow 0} H(s) = 0
\quad \text{and} \quad
\lim _{s \uparrow \infty} H(s) = H_\infty .
\nonumber
\end{gather}

\end{cor}
\noindent {\bf Proof.}
The result follows immediately from
Theorem~\ref{prop:Q}.(III) and the fact that the
right-hand side of (\ref{H-ODE}) is strictly
positive for all values of $H(s)$ in
$\bigl] 0 , G(s) \bigr[$.
In particular, the definition of $Q$ implies that
\be
\lim _{s \downarrow 0} H(s) =
\lim _{s \downarrow 0} F(s) Q(s) = 0
\quad \text{and} \quad
\lim _{s \uparrow \infty} H(s) =
\lim _{s \uparrow \infty} F(s) Q(s) = H_\infty ,
\ee
for $H_\infty = q_\infty$, where $q_\infty \in
\mbox{} ]0 , G_\infty]$ is as in
Theorem~\ref{prop:Q}.(III).
\mbox{}\hfill$\Box$

\section{The solution to the optimal stopping problem}
\label{sec:w}

Each of the solutions to the ODE (\ref{H-ODE})
derived in Corollary~\ref{cor:H-ODE} is associated
with a function $w$ given by (\ref{w}) that is a
candidate for the optimal stopping problem's value
function $v$.
The following result presents a comprehensive
study of these functions $w$.
It turns out that the point
\ben
H_\infty^\circ = \frac{m+1}{m} \in \mbox{}
]0 , G_\infty[ \label{Hoinfty}
\een
identifies the free-boundary function that yields
the solution to the optimal stopping problem.

\begin{thm} \label{thm:w-VI}
Suppose that Assumption~\ref{A} holds true.
Also, consider the function $w$ given by
(\ref{w}) with $A$ and $B$ given by (\ref{AnB})
and for $H$ being any of the solutions to the
ODE (\ref{H-ODE}) that are as in
Corollary~\ref{cor:H-ODE}.
The following statements hold true.
\smallskip

\noindent {\rm (I)}
The function $w$ is $C^1$ and its restriction in
\be
\bigl\{ (x,s) \in \bbr_+^2 \mid \ 0 < x \leq s
\text{ and } x \neq H(s) \bigr\}
\ee
is $C^2$.
\smallskip

\noindent {\rm (II)}
If $H_\infty \in \mbox{} ]H_\infty^\circ, G_\infty]$,
then $w$ does not satisfy the variational inequality
(\ref{HJB}) because there exist $0 < x \leq s$
such that $w(x,s) < 0$.
\smallskip

\noindent {\rm (III)}
If $H_\infty \in \mbox{} ]0, H_\infty^\circ]$,
then $w$ is strictly positive and satisfies the
variational inequality (\ref{HJB}) as well as the
boundary condition (\ref{HJB-BC}).
\smallskip

\noindent {\rm (IV)}
If $H_\infty \in \mbox{} ]0, H_\infty^\circ[$,
then $w$ does not satisfy the transversality
condition (\ref{TV-cond}).
Moreover, if $(\tau_\ell)$ is any sequence
of bounded $(\fcal _t)$-stopping times such
that $\lim _{\ell \uparrow \infty} \tau_\ell
= \infty$, then
\be
\lim _{\ell \uparrow \infty} \EXP \bigl[
e^{-r \tau_\ell} w(X_{\tau_\ell} , S_{\tau_\ell})
\bigr] \geq A_\infty x^n > 0 ,
\ee
for some constant $A_\infty = A_\infty (H_\infty)
> 0$.
\smallskip

\noindent {\rm (V)}
If $H_\infty = H_\infty^\circ$, then $w$ satisfies
the transversality condition (\ref{TV-cond}).
Furthermore, if $(\tau_\ell)$ is any sequence
of bounded $(\fcal _t)$-stopping times such
that $\lim _{\ell \uparrow \infty} \tau_\ell
= \infty$, then
\ben
\lim _{\ell \uparrow \infty} \EXP \bigl[
e^{-r \tau_\ell} w(X_{\tau_\ell} , S_{\tau_\ell})
\bigr] = 0 . \label{TV-cond-S}
\een
\end{thm}
\noindent {\bf Proof.}
The claims in (I) follow from the construction of
$w$ (see also the first paragraph of Section~\ref{sec:H}).
\smallskip

{\em Proof of part\/} (II).
Differentiating the expression for $A$ given by
(\ref{AnB}) and using the ODE (\ref{H-ODE}), we
obtain
\ben
\dot{A} (s) = - \frac{\dot{F} (s) \bigl( H(s)/s \bigr)^{n-m}}
{H^{n+1} (s) \Bigl( 1 - \bigl( H(s)/s \bigr)^{n-m} \Bigr)}
< 0 \quad \text{for all } s > 0 . \label{A'}
\een
On the other hand, passing to the limit as $s \uparrow
\infty$ in the same expression yields
\ben
A_\infty := \lim _{s \uparrow \infty} A(s) =
\frac{m H_\infty - (m+1)}{(n-m) H^{n+1}_\infty} 
\begin{cases}
> 0, & \text{if } H_\infty \in \mbox{} ]0, H_\infty^\circ[ , \\
= 0, & \text{if } H_\infty = H_\infty^\circ , \\
< 0, & \text{if } H_\infty \in \mbox{} ]H_\infty^\circ ,
G_\infty] . \end{cases} \label{Ainfty}
\een
In view of these calculations, we can see that
\begin{align}
\text{if } \lim _{s \uparrow \infty} H(s) = H_\infty
\in \mbox{} ]0, H_\infty^\circ] , \quad & \text{then }
A(s) > 0 \quad \text{for all } s > 0 \label{A>0} \\
\text{and, if } \lim _{s \uparrow \infty} H(s) = H_\infty
\in \mbox{} ]H_\infty^\circ, G_\infty] , \quad & \text{then }
A(s) \begin{cases} > 0 , & \text{for all } s \in \mbox{}
\bigl] 0,\overline{s}[, \\ < 0 , & \text{for all } s >
\overline{s} , \end{cases} \label{A<0}
\end{align}
for some $\overline{s} = \overline{s} (H_\infty)
\geq 0$.
Similarly, we calculate
\begin{gather}
\dot{B} (s) = \frac{\dot{F} (s)}
{H^{m+1}(s) \Bigl( 1 - \bigl( H(s)/s \bigr)^{n-m} \Bigr)}
> 0 \quad \text{for all } s > 0 \nonumber \\
\text{and} \quad
0 < B(s) < B_\infty := \lim _{s \uparrow \infty} B(s)
= \frac{n+1 - n H_\infty}{(n-m) H_\infty^{m+1}} 
\quad \text{for all } s>0 . \label{Binfty}
\end{gather}
The claims in (\ref{A<0}) and (\ref{Binfty})
imply that
\be
\text{if } \lim _{s \uparrow \infty} H(s) = H_\infty
\in \mbox{} ]H_\infty^\circ, G_\infty] , \quad \text{then }
\lim _{s \uparrow \infty} \frac{B(s)}{A(s)} \in \mbox{}
] {-\infty}, 0[ .
\ee
Part~(II) of the theorem follows from this observation
and the fact that, given any $s > \overline{s}$,
\be
w(x,s) < 0 \quad \Leftrightarrow \quad
x^{n-m} > - \frac{B(s)}{A(s)} ,
\ee
where $\overline{s}$ is as in (\ref{A<0}).
\smallskip

{\em Proof of part\/} (III).
Suppose that the free-boundary function $H$ is such
that $H_\infty \in \mbox{} ]0, H_\infty^\circ]$.
In this case, (\ref{A>0}) and (\ref{Binfty}) imply that
$w$ is strictly positive.
In view of this observation and its construction, we
will prove that $w$ satisfies the variational inequality
(\ref{HJB}) with boundary condition (\ref{HJB-BC})
if we show that
\begin{align}
f(x,s) & := \frac{1}{2} \sigma^2 x^2
\frac{\partial ^2 \bigl( x^{-1} F(s) - 1 \bigr)}{\partial x^2}
+ \mu x \frac{\partial \bigl( x^{-1} F(s) - 1 \bigr)}
{\partial x} - r \bigl( F(s) x^{-1} - 1 \bigr) \nonumber \\
& = (\sigma^2 - \mu - r) x^{-1} F(s) + r
\leq 0 \quad \text{for all } s > 0 \text{ and } x \in
\bigl] 0, H(s) \bigr[ \label{HJB-S} \\
g(x,s) & := w(x,s) - x^{-1} F(s) + 1 \geq 0 
\quad \text{for all } s > 0 \text{ and } x \in
\bigl] H(s) , s \bigr[. \label{HJB-W}
\end{align}
The inequality (\ref{HJB-S}) follows immediately
from (\ref{G}) and the fact that $0 < H(s) < G(s)$
for all $s>0$.

The strict positivity of $w$ implies that (\ref{HJB-W})
holds true for all $s>0$ and $x \in \mbox{} \bigl[
F(s) , s \bigr[$.
To show that the inequality holds true for all
$s>0$ and $x \in \mbox{} \bigl] H(s) , F(s)
\bigr[$, we first note that
\be
\frac{1}{2} \sigma^2 x^2 g_{xx} (x,s) + \mu x
g_x (x,s) - rg(x,s) = - f(x,s) \quad \text{for all }
s > 0 \text{ and } x \in \mbox{} ]H(s), s[ ,
\ee
where $f$ is defined as in (\ref{HJB-S}).
Combining this identity with the inequalities
\be
- f(x,s) = \begin{cases} 
> 0 , & \text{if } x \in \mbox{} ]H(s), G(s)[  , \\ 
< 0 , & \text{if } x \in \mbox{} ]G(s), s[ , \end{cases}
\ee
which follow from (\ref{G}), (\ref{G(s)<s}), the fact
that $0 < H(s) < G(s)$ for all $s>0$ (see
Corollary~\ref{cor:H-ODE}) and the maximum
principle, we can see that, given any $s > 0$,
\begin{align}
\text{the function } g(\cdot , s) \text{ has }
\begin{cases} 
\text{no positive maximum inside } ]H(s), G(s)[ , \\
\text{no negative minimum inside } ]G(s), s[ . 
\end{cases} \label{MP1}
\end{align}

In view of the limit 
\begin{align}
\lim _{x \downarrow H(s)} g_{xx} (x,s) = \mbox{}
& n (n-1) \frac{m H(s) - (m+1) F(s)}{n-m} H^{-3}
(s) \nonumber \\
& + m (m-1) \frac{(n+1) F(s) - n H(s)}{n-m} H^{-3}
(s) - 2 F(s) H^{-3} (s) \nonumber \\
= & - mn H^{-3} (s) \bigl( G(s) - H(s) \bigr) > 0 , \nonumber
\end{align}
where the inequality follows from Corollary~\ref{cor:H-ODE},
and the identities 
$g \bigl( H(s) , s \bigr) = g_x \bigl( H(s), s \bigr) = 0$,
which follow from the $C^1$-continuity of $w(\cdot , s)$
along $H(s)$, we can see that
\be
g_x \bigl( H(s) {+} \varepsilon , s \bigr) > 0 
\quad \text{and} \quad 
g \bigl( H(s) {+} \varepsilon , s \bigr) > 0 
\quad \text{for all } \varepsilon > 0
\text{ sufficiently small} .
\ee
Combining this observation with (\ref{MP1})
and the fact that $g \bigl( F(s), s \bigr) > 0$,
we obtain (\ref{HJB-W}) for all $s>0$ and
$x \in \mbox{} \bigl] H(s) , s \bigr[$.
\smallskip

{\em Proof of parts\/} (IV) {\em and\/} (V):
{\em Preliminary analysis\/}.
Let $(\tau_\ell)$ be a sequence of bounded
$(\fcal _t)$-stopping times such that $\lim
_{\ell \uparrow \infty} \tau_\ell = \infty$.
In view of the observation that
\be
0 \leq w(x,s) {\bf 1} _{\{ x \leq H(s) \}} = \bigl( x^{-1}
F(s) - 1 \bigr) {\bf 1} _{\{ x \leq H(s) \}} \, \leq
x^{-1} 
\ee
and (\ref{TV-F}), we can see that
\ben
\lim _{\ell \uparrow \infty} \EXP \bigl[
e^{-r \tau_\ell} w(X_{\tau_\ell} , S_{\tau_\ell})
\bigr] = \lim _{\ell \uparrow \infty} \EXP \bigl[
e^{-r \tau_\ell} w(X_{\tau_\ell} , S_{\tau_\ell})
{\bf 1} _{\{ X_{\tau_\ell} > H(S_{\tau_\ell}) \}}
\bigr] . \label{TV-P1}
\een

Recalling Assumption~\ref{A} and using the
expression for $A$ given by (\ref{AnB}) as well
as the definition (\ref{Hoinfty}) of
$H_\infty^\circ$, we obtain
\ben
0 \leq A(s) x^n  {\bf 1}_{\{ x > H(s) \}} <
\frac{m H(s) - (m+1)}{(n-m) H^{n+1}(s)} s^n
= \frac{-m}{(n-m) H^{n+1}(s)} s^n \bigl(
H_\infty^\circ - H(s) \bigr) . \label{Ax^n}
\een
On the other hand, (\ref{Binfty}) implies that
\ben
0 \leq B(u) x^m  {\bf 1}_{\{ x > H(u) \}} < B_\infty
H^m (s) \quad \text{for all } u \geq s , \label{Bx^m}
\een
where $B_\infty$ is defined as in (\ref{Binfty}).
\smallskip

{\em Proof of part\/} (IV).
Fix any $H_\infty \in \mbox{} ]0, H_\infty^\circ[$
and let $H$ be the solution to the ODE
(\ref{H-ODE}) satisfying $\lim _{s \uparrow \infty}
H(s) = H_\infty$.
In view of (\ref{A'}), (\ref{Ainfty}), the strict positivity
of $B$, (\ref{TV-P1}) and the inclusion $\{ X_{\tau_\ell}
> H_\infty \} \subseteq \bigl\{ X_{\tau_\ell} >
H(S_{\tau_\ell}) \bigr\}$, we can see that
\be
\lim _{\ell \uparrow \infty} \EXP \bigl[ e^{-r \tau_\ell}
w(X_{\tau_\ell} , S_{\tau_\ell}) \bigr] \geq A_\infty
\lim _{\ell \uparrow \infty} \EXP \bigl[ e^{-r \tau_\ell}
X_{\tau_\ell}^n {\bf 1} _{\{ X_{\tau_\ell} > H_\infty \}}
\bigr] .
\ee
Combining this observation with the fact that
\be
0 \leq \EXP \bigl[ e^{-r \tau_\ell} X_{\tau_\ell}^n
{\bf 1} _{\{ X_{\tau_\ell} \leq H_\infty \}} \bigr]
\leq H_\infty^n \EXP \bigl[ e^{-r \tau_\ell} \bigr]
\xrightarrow[\ell \uparrow \infty]{} 0 ,
\ee
we obtain
\be
\lim _{\ell \uparrow \infty} \EXP \bigl[
e^{-r \tau_\ell} w(X_{\tau_\ell} , S_{\tau_\ell})
\bigr] \geq A_\infty \lim _{\ell \uparrow \infty}
\EXP \bigl[ e^{-r \tau_\ell} X_{\tau_\ell}^n
\bigr] = A_\infty x^n > 0 .
\ee

{\em Proof of part\/} (V).
Let $H^\circ$ be the solution to the ODE
(\ref{H-ODE}) satisfying $\lim _{s \uparrow \infty}
H^\circ(s) = H_\infty^\circ$.
Using L'H\^{o}pital's lemma and the definition
(\ref{G}) of $G$, we calculate
\begin{align}
\lim_{s \uparrow \infty} s^n \bigl( H_\infty^\circ
- H^\circ (s) \bigr) & = \lim_{s \uparrow \infty}
\frac{\dot{H}^\circ (s)}{ns^{-n-1}} \nonumber \\
& = \lim_{s \uparrow \infty}
\frac{\dot{F} (s)}{ns^{-n-1}}
\frac{\Bigl( (n+1) \bigl( H^\circ(s) / s \bigr)^{n-m}
- (m+1) \Bigr) H^\circ(s)}
{-mn \bigl( G(s) - H^\circ (s) \bigr) \Bigl( 1 - \bigl(
H^\circ (s) / s \bigr) ^{n-m} \Bigr)} = 0 \nonumber
\end{align}
because
\be
\lim_{s \uparrow \infty}
\frac{\dot{F} (s)}{ns^{-n-1}}
= \frac{1}{n} \lim_{s \uparrow \infty} s^{n+1} e^{-s} = 0 . 
\ee
Therefore,
\be
\max _{u \geq s} u^n \bigl( H_\infty^\circ
- H^\circ (u) \bigr) < \infty .
\ee
Combining this observation with (\ref{Ax^n})
and (\ref{Bx^m}), we can see that
\be
\max _{x>0 , \, u \geq s} w(x,u) {\bf 1}
_{\{ H^\circ (u) < x \leq u \}} < \infty .
\ee
The claims in part~(V) of the theorem follow
from this result and (\ref{TV-P1}). 
\mbox{}\hfill$\Box$
\bigskip

The following result provides the solution to the
optimal stopping problem considered in this
paper.

\begin{thm} \label{thm:main}
Consider the optimal stopping problem defined by (\ref{v-May}), (\ref{X})--(\ref{F}) 
and suppose that the problem's
data satisfy Assumption~\ref{A}.
The problem's value function $v$ identifies with the
function $w$ defined by (\ref{w}) for $H=H^\circ$
being the solution to the ODE (\ref{H-ODE})
characterised by $\lim _{s \uparrow \infty} H^\circ
(s) = H_\infty^\circ$, while
\ben
\tau_\star = \inf \bigl\{ t \geq 0 \mid \ (X_t , S_t) \in
\sscr \bigr\} = \inf \bigl\{ t \geq 0 \mid \ X_t \leq
H^\circ (S_t) \bigr\} \label{opt-tau}
\een
is an optimal stopping time.
\end{thm}
\noindent {\bf Proof.}
Fix any $s>0$ and $x \in \mbox{} ]0,s]$.
Using It\^{o}'s formula and the fact that $S$ increases
inside the set $\{ X = S \}$, we obtain
\begin{align}
e^{-rT} w (X_T,S_T) = \mbox{} & w (x,s) + \int _0^T e^{-rt}
w_s (S_t,S_t) \, \di S_t \nonumber \\
& + \int _0^T e^{-rt} \biggl( \frac{1}{2} \sigma^2 X_t^2 w_{xx}
(X_t,S_t) + \mu X_t w_x (X_t,S_t) - r w (X_t,S_t) \biggr) \, \di t
+ M_T , \nonumber
\end{align}
where
\be
M_T = \sigma \int _0^T e^{-rt} X_t w_x (X_t,S_t) \, \di W_t .
\ee
Therefore,
\begin{align}
e^{-rT} \bigl( X_T^{-1} & F(S_T) - 1 \bigr)
\nonumber \\
= \mbox{} & w (x,s) + e^{-rT} \bigl( X_T^{-1}
F(S_T) - 1 - w (X_T,S_T) \bigr) + \int _0^T
e^{-rt} w_s (S_t,S_t) \, \di S_t \nonumber \\
& + \int _0^T e^{-rt} \biggl( \frac{1}{2} \sigma^2
X_t^2 w_{xx} (X_t,S_t) + \mu X_t w_x (X_t,S_t)
- r w (X_t,S_t) \biggr) \, \di t + M_T . \label{ITO}
\end{align}
Given a stopping time $\tau \in \tcal$ and a localising
sequence of bounded stopping times $(\tau_\ell)$
for the local martingale $M$, this identity and
the fact that $w$ satisfies the variational inequality
(\ref{HJB}) as well as the boundary condition
(\ref{HJB-BC}) imply that
\ben
\EXP \Bigl[ e^{-r (\tau \wedge \tau_\ell)}
\bigl( X_{\tau \wedge \tau_\ell}^{-1}
F(S_{\tau \wedge \tau_\ell}) - 1 \bigr) \Bigr]
\leq w (x,s) . \label{VT00}
\een
Furthermore, Fatou's lemma implies that
\ben
\EXP \Bigl[ e^{-r \tau} \bigl( X_{\tau}^{-1}
F(S_{\tau}) - 1  \bigr) {\bf 1} _{\{ \tau < \infty \}} 
\Bigr] \leq \liminf _{\ell \uparrow \infty}
\EXP \Bigl[ e^{-r (\tau \wedge \tau_\ell)} \bigl(
 X_{\tau \wedge \tau_\ell}^{-1}
 F(S_{\tau \wedge \tau_\ell}) - 1 \bigr) \Bigr]
 \leq w (x,s) . \label{VT-V<w}
\een
It follows that $v(x,s) \leq w(x,s)$.

If $\tau_\star \in \tcal$ is the stopping time defined
by (\ref{opt-tau}), then (\ref{ITO}) implies that
\be
\EXP \Bigl[ e^{-r \tau_\star} \bigl( X_{\tau_\star}^{-1}
F(S_{\tau_\star}) - 1 \bigr) {\bf 1}
_{\{ \tau_\star \leq \tau_\ell \}} \Bigr] = w (x,s)
- \EXP \Bigl[ e^{-r \tau_\ell} w (X_{\tau_\ell},
S_{\tau_\ell}) {\bf 1} _{\{ \tau_\ell < \tau_\star \}}
\Bigr] . 
\ee
Passing to the limit as $\ell \uparrow \infty$ using
the monotone convergence theorem, the positivity
of $w$ and (\ref{TV-cond-S}), we obtain 
\be
\EXP \Bigl[ e^{-r \tau_\star} \bigl( X_{\tau_\star}^{-1}
F(S_{\tau_\star}) - 1\bigr) {\bf 1}
_{\{ \tau_\star < \infty \}} \Bigr] = w (x,s) ,
\ee
which implies that $v(x,s) \geq w(x,s)$.
This result and the inequality $v(x,s) \leq
w(x,s)$, which follows from (\ref{VT-V<w}), imply the
claims of the theorem.
\mbox{}\hfill$\Box$

\section*{Acknowledgment}

We are grateful to Goran Peskir for several insightful
and valuable discussions.

\small

\end{document}